\documentclass{scrartcl}
\usepackage[T1]{fontenc}
\usepackage[latin1]{inputenc}
\usepackage{enumerate}
\usepackage{a4,amssymb,amsmath}
\usepackage[mathscr]{eucal}

\newtheorem{thm}{Theorem}[section]

\newtheorem{lem}[thm]{Lemma}

\newtheorem{defn}[thm]{Definition}

\newcommand{\norm}[1]{\left\Vert#1\right\Vert}

\newcommand{\Cp}{\mathbb{C}}

\renewcommand{\rho}{\varrho}
\renewcommand{\phi}{\varphi}

\newcommand{\C}{\mathcal{C}}

\newcommand{\med}{\medskip \\}
\newcommand{\bg}{\bigskip \\}

\newcommand{\be}{\begin{equation}}
\newcommand{\ee}{\end{equation}}
\newcommand{\bes}{\begin{equation*}}
\newcommand{\ees}{\end{equation*}}

\newcommand{\tnorm}[1]{\mathrel{%
 \left\vert\kern-1pt\left\vert\kern-1pt\left\vert%
  #1%
 \right\vert\kern-1pt\right\vert\kern-1pt\right\vert}}

\newcommand{\tnorma}[1]{\mathrel{%
 \vert\kern-1pt\vert\kern-1pt\vert%
  #1%
 \vert\kern-1pt\vert\kern-1pt\vert}}

\newcommand{\tnormb}[2]{\mathrel{%
 #2\vert\kern-1pt#2\vert\kern-1pt#2\vert%
  #1%
 #2\vert\kern-1pt#2\vert\kern-1pt#2\vert}}

\def\carre{\hbox{
\vrule height 1.453ex  width 0.093ex  depth 0ex
\vrule height 1.5ex  width 1.3ex  depth -1.407ex\kern-0.1ex
\vrule height 1.453ex  width 0.093ex  depth 0ex\kern-1.35ex
\vrule height 0.093ex  width 1.3ex  depth 0ex}\,}
\def\buildo#1^#2{\mathrel{\mathop{\null#1}\limits^{#2}}}
\def\buildu#1_#2{\mathrel{\mathop{\null#1}\limits_{#2}}}

\begin{document}

\title{Boundary Behavior of the Bergman Metric}%
\author{Klas Diederich and J.E.Fornaess}%


\begin{abstract}
We prove more precise estimates for the Bergman metric on strongly pseudoconvex domains, based on the squeezing function of the domain.
\end{abstract}
\maketitle
\section{Introduction}

The goal of this paper is to prove more precise estimates on the boundary behavior of the Bergman metric than known so-far. We consider the Bergman metric on a bounded strictly pseudoconvex domain ${\Omega }$ with $\C^4$ boundary. (Our result becomes weaker, if we even allow $\C^3$ boundaries) However, below this degree of differentiability we cannot say anything with our methods.

We thank N. Nikolov for pointing out to us the similarities between his paper in J.Math.Anal.Appl. 421 (2015) 180-185 and our announcement.

The major tool is the consideration of the squeezing function. It is defined as follows: for a given injective holomorphic map $f:\Omega \rightarrow \textbf{B} {^n}$ satisfying $f(z)=0$ we set

$$ S_{\Omega,f}:=sup \{r > 0: r\mathbb{B}^n\subset f(\Omega)\} $$
and we put
$$ S_\Omega (z) := \sup_f \{S_{\Omega,f}(z)\} $$
where $f$ ranges over all injective holomorphic maps $f:\Omega \rightarrow \mathbb{B}^n$ with $f(z)=0$
if $\Omega$ is a bounded domain in $\mathbb C^n$.
In [6] the following Theorem has been proved:
\begin{thm}\label{T:M} Let $\Omega = \{\delta < 0 \}\subset \Cp ^n$ be a strictly pseudoconvex domain with a defining function $\delta$ of class $\mathcal{C} ^k$ for $k\geq 3$. The squeezing function, see Deng--Guan--Zhang [1],  $S_\Omega (z)$ for $\Omega$ satisfies the estimate
$$S_{\Omega} (z) \geq 1-C \sqrt {|\delta (z)|} $$
for a fixed constant $C$. If we even have $k\geq 4$, then there exists a constant $C>0$ such that the squeezing function $S_\Omega (z)$ for $\Omega$ satisfies
$$  S_\Omega (z)\geq 1-C|\delta (z)|$$
for all $z$.
\end{thm}

In [6] this was used to give boundary estimates for the Caratheodory metric. This was based on earlier work by Ma [9], see also Fu [8] and Diederich--Forn\ae ss--Wold [5]. We will in this paper use this result to provide similar estimates for the Bergman metric.
See also Diederich [2], Diederich [3], Diederich--Forn\ae ss-Herbort [4], Fu[7].

Our main theorem is the following:
\begin{thm}\label{metrics}
Let $\Omega\subset\mathbb C^n$ be a strictly pseudoconvex domain of class $\mathcal C^3$,
let $p\in b\Omega$, and let $\delta$ be a defining function for $\Omega$ near $p$,
such that $\|\nabla\delta(z)\|=1$ for all $z\in b\Omega$.  Then if $d_\Omega(z,\zeta)$
is the Bergman metric, there exists a constant $C>0$
such that
\begin{align*}
&(n+1) (1-C\sqrt{|\delta(z)|})\left[\frac{\mathfrak{L}_{\pi(z)}(\xi_T)}{|\delta(z)|} + \frac{\|\xi_N\|}{4\delta(z)^2}\right]^{1/2}\leq d_\Omega(z,\xi)\\
& \leq (n+1) (1+C\sqrt{|\delta(z)|})\left[\frac{\mathfrak{L}_{\pi(z)}(\xi_T)}{|\delta(z)|} + \frac{\|\xi_N\|}{4\delta(z)^2}\right]^{1/2}
\end{align*}
for all $z$ near $p$, and all $\xi=\xi_N+\xi_T$, where $\pi$ is the orthogonal projection
to $b\Omega$, $\xi_N$ is the complex normal component of $\xi$ at $\pi(z)$ and $\xi_T$ is the complex tangential component, and $\mathfrak{L}$ is the Levi form of $\delta$.
\end{thm}

The theorem gives an estimate for the Bergman metric which is valid on any strongly pseudoconvex domain with $\mathcal{C}^3$ boundary.
In the case of $\mathcal{C}^4$ boundary, it is possible to improve the estimates if one first changes coordinates appropriately.
We can then prove, using the coordinates explained in Section 4.

\begin{thm}[Sharp estimates for the Bergman metric]
Put $B:= B{_\Omega} ((r,0),\xi)$, where $0\in b\Omega$ and $r>0$ is a small positive radius on the inner normal to $b\Omega$ at $0$. Then we have, if the boundary is $\mathcal
{C}^4,$
$$B\leq (n+1) (1+Cd) \sqrt{\frac{|\xi_N|^2}{4d^2}+\frac{L(\xi_T)}{d}}$$
In the opposite sense we have the estimates
\begin{equation}
B \geq
\begin{cases}
(n+1)(1-Cd)\sqrt{{\frac{|\xi_N|^2}{4d^2}+\frac{L(\xi_T)}{d}}} &\text{if the boundary regularity is $\mathcal{C}^4 $}\\
\end{cases}
\end{equation}
\end{thm}

We observe  that these estimates are then also valid for the Kobayashi, Sibony, Azukawa and Caratheodory metric.

\section{Basics from the Bergman Theory}
First of all we remind the reader of the following monotonicity result for the Bergman kernel function:

\begin{lem}
Let $0\in \mathbb{B}(0;r_1) \subset \Omega \subset \mathbb{B}(0;r_2). $
Then we have for the Bergman kernel function:
$$K_{\mathbb{B}(0,r_1)}(0) \geq K_\Omega (0)\geq K_{\mathbb{B}(0,r_2)}(0) $$
\end{lem}
It is typical for the Bergman theory, that certain Maxima of linear evaluation functionals play an important role. We define
\begin{defn}
For a bounded domain $\Omega$, point $0\in \Omega$ and holomorphic tangent vectors $\xi$ at $0$ to $\Omega$ we define
$$M_{\Omega} (0,\xi):= \max_{\norm f _{L^2} =1} \norm{ f'(\xi)}(0) $$
\end{defn}
It is easily seen, that the minimum functional as defined in the last definition also satisfies a monotonicity property similar to the Bergman kernel itself as expressed in the following Lemma:
\begin{lem}
Let $\Omega $ be a domain with $\mathbb{B}(0,r_2) \supseteq \Omega \supseteq \mathbb{B}(0,r_1)\ni $, then one has
$$M_{\mathbb{B}(0,r_1)}(0,\xi)\geq M_\Omega (0,\xi)\geq M_{\mathbb{B}(0,r_2)}(0,\xi)$$
\end{lem}
We still need to mention the explicit forms of the Bergman kernel of the ball of radius R and also the respective Maximum as defined in Definition 3.2. Both can easily be calculated. We get:

For the ball of radius $r$ we have the following explicit formula for the Bergman kernel evaluated at $0$:
\begin{lem}
$K_{\mathbb{B}(0,r)} = \frac{n!}{\pi ^n r^{2n}}$
\end{lem}

\begin{lem}
Let now $d_\Omega (p,\xi)$ be the Bergman metric of a domain $\Omega$  at the point $p$ applied to the vector $\xi.$
Furthermore, let $M(p,\xi)$ be the linear functional $M$ defined in Def. 3.2 at the point $p$ and the vector $\xi$. Then we have the following well-known formula:
$$d_\Omega(p,\xi)= \frac{M(p,\xi}{\sqrt{(K(p))}}$$
where K means the Bergman kernel function at the point $p$
\end{lem}

Furthermore, it is easy to see that we have the following explicit formula for the $M$ from Def. 3.2 evaluated at $0$ and a vector $\xi$:
\begin{lem}
$ M_{\mathbb B(0,r)}(0,\xi)=(n+1) \sqrt{{\frac{n!}{\pi ^n}}\frac{1}{r^{n+2}}} \norm{\xi} ^2$
\end{lem}
Finally, we get the following explicit estimate
\begin{lem}\label{L:ME}
$$\frac{M_{B}(0,r_1)}{\sqrt{K_{B(0,r_2)}}}\geq d_\Omega (0,\xi) = \frac{M_\Omega}{\sqrt{K_\Omega}}\geq \frac{M_\mathbb{B}(0,r_2)}{\sqrt{K_\mathbb{b}(0,r_1)}}  $$
\end{lem}

\section{Proof of Theorem 1.2}

\begin{thm}
Let $\Omega$ be a bounded domain in $\mathbb C^n.$ Let $d^K_\Omega$ denote the Kobayashi metric.
We then have
$$
(n+1)d_\Omega^K(p,\xi) S_\Omega(p)^{\frac{n+2}{2}}
\leq d_\Omega(p,\xi)\leq (n+1)d_\Omega^K(p,\xi) \frac{1}{S_\Omega(p)^{\frac{n+2}{2}}}
$$
\end{thm}

To prove the Theorem, we observe that both the Bergman metric and Kobayashi metric are invariant under biholomorphic maps. Hence it suffices to consider the case when $\Omega $ is contained on the unit ball,
and contains the ball of radius $S_\Omega(z)$ and we estimate the metrics at the origin.
But then we can use Lemmas 2.4, 2.6 and 2.7.

Now Theorem 1.2 follows from the result by Ma for the Kobayashi metric.

We remark that the estimate in Theorem 1.1 for the squeezing function iin the case $\mathcal {C}^3$ is all we need.
 However, in the $\mathcal {C}^4$ case we get a sharper comparison with the Kobayashi metric
similar to the ones for the Caratheodory metric, Azukawa metric and Sibony metric in [6].

\section{Exposing points and the squeezing function}
We need a more precise geometric setup. Let $k=3$ or $4$ and let $\Omega$ be a bounded strongly pseudoconvex domain of class $\mathcal{C}^k$. In Lemma 5.1 of [6] the following result is proved:

\begin{lem}
Let $p\in b\Omega$. There exists an injective holomorphic map $\Phi : \bar\Phi\rightarrow \Cp ^n$, such that $\tilde{\Omega} = \Phi (\Omega)$ satisfying the following:
\begin{itemize}
\item[(i)] $\tilde\Omega\subset\mathbb B^n$,
\item[(ii)] $\phi(p)=(1,0,\cdot\cdot\cdot,0)=:a$ and $\phi^{-1}(b\mathbb B^n)=\{p\}$,
\item[(iii)] near $a$ we have that, $\tilde\Omega=\{\rho<\mu^2\}, 0<\mu<1$ where
$$
\rho(z)=|z_1-(1-\mu)|^2+\|z'\|^2 +  O(|z_1-1|^2) + O(\|z-a\|^k).
$$
\end{itemize}
\end{lem}
\underline{Proof.} The details of the proof can be found in [6].We do not repeat them here.\med
A major step in the construction of the squeezing function consists in osculating the given strictly pseudoconvex domain by the ball in the following way:

\begin{lem}
There exists a constant $C>0$ such that for $w\in b\Phi_r (\Omega)$ we have
$$|w_1|^2+\frac{1}{\mu} \norm {w'}^2 \leq 1+C(1-r)^\frac{k-2}{2}$$
\end{lem}
Furthermore, Lemma 5.4 in [6] tells us:

 \begin{lem}
We set
$\tilde{\eta} = \begin{cases}
		\eta, & k=4 \\
		\frac{\eta}{1-C\eta},& k=3
	\end{cases}$\med
(i) If $k=4$ then $B^{\mu}_{\eta,\tilde{\eta}}\subset \Omega$ for all $\eta$ small enough\\
(ii) If $k=3$,
and the constant $C>0$ is fixed large enough, then $B^{\mu}_{\eta,\tilde{\eta}}\subset \Omega$ for all $\eta$ small enough.
\end{lem}
In Lemma 5.6 of [6] it was furthermore proved using the previous Lemmas:
\begin{lem}
Let $\psi (z)=(z_1, \frac{1}{\sqrt{\mu}z'}$. Suppose that $0<\eta, r>1, 1-2\eta < r$ and $\tilde{\eta}>0$. Then
$\psi(\phi_r (B^\mu _{\eta,\tilde{\eta}}))$ contains the ball of radius
$$\sqrt{1-2(1-r)\frac{1}{\tilde\eta }-4\cdot |1-\frac{\eta}{\tilde{\eta}}}| $$
\end{lem}

\section{Proof of Theorem 1.3}

We choose a point of the form $(r,0)$ near the boundary of $\Omega$. We decompose the tangent vector $\xi$ to $\Omega$ at $(r,0)$ in the form
$$\xi = \xi_T + \xi_N $$
where T stands for the tangential directions and N is the normal direction. We consider the mapping
$$\phi_r (z_1,z'):= ( \frac{z_1 -r}{1-z_1 r} ,\frac{\sqrt{1-r^2}}{1-z_1 r} z' ) $$

After composing this map with a stretching factor $\frac{1}{\sqrt{\mu}}$, in the $\xi'$ variable, we obtain the following map:
$$\psi(\phi_r)= (\frac{z_1-r}{1-z_1 r}, \frac{1}{\sqrt{\mu}}\frac{\sqrt{1-r^2}}{1-z_1 r}  z') $$
It is the final form of our adaptation to the use of the squeezing map.\med
Calculating this map at the point $(r, \xi)$ we get
$$\lambda = \Bigl( \frac{1}{1-r^2 }\xi _1 , \frac{1}{\sqrt{\mu}} \frac{1}{1-r^2}\xi ' \Bigr) $$
From this we get the following splitting into normal and tangential component:
$$\| \lambda \| = \sqrt{\frac{1}{(1-r^2)}^2}\|\xi _1|^2 + \frac{1}{\mu} \frac{1}{1-r^2} \| \xi '\|^2 $$
The point $r$ has boundary distance $d=1-r$. Therefore we have $1-r^2=(1+r)(1-r)$. Putting this into our expression for $\lambda$, we get:
$$\|\lambda \|={\sqrt{\frac{|\xi_N| ^2}{d^2}(\frac{1}{1+r ^2} + \frac{1}{\mu}\frac{1}{d}\frac{1}{1+r} \|\xi _{T} \| ^2 }}) $$
\med
Our next goal is, to write a defining function $\rho$ for the boundary of $\Omega$ in terms of a normalized defining function. We do this in the following way:
$\Omega$ can be written as
$$\Omega=\{\rho - \mu^2<0\}$$
More precisely, $\rho$ has the form
$$\rho = |z_1-1|^2 + 2 Re(\mu (z_1-1))+\|z'\|^2+.... $$
This can be written as:
$$\rho=\|z\|^2+ z\mu x_1 + ... $$
with $| \nabla \rho | = 2\mu $. From this we easily get:
$$\mathfrak{L}_{\rho/(2\mu})(\xi _T )= \frac{1}{2\mu}\|\xi_T \| ^2   $$

We now express $\lambda$ introducing the Levi form into it and obtain:
$$\|\lambda\| = ({ \sqrt{{\frac{|\xi _N|^2}{d^2}}\frac{1}{(2-d)^2} + \frac{2}{d}\frac{1}{2-d} \mathfrak{L}_{\rho/(2\mu}(\xi _T )}}  $$
In other form
$$\|\lambda\| = {\sqrt{\frac{1}{4}(\frac{1}{1-d/2})^2\frac{|\xi _N|^2}{d^2} + \frac{1}{d}\frac{1}{1-d/2} \mathfrak{L}(\xi _T)}}  $$
This finally gives
$$\|\lambda\|=(1\pm d)\sqrt{\frac{|\xi _N|^2}{4d{^2}} + \frac{L(\xi_T)}{d}}$$
\bg
In order to make further progress in estimating the behavior of the Bergman metric, we next use extensively the inequalities stated in Lemma \ref{L:ME}. They can now be stated explicitly for $r_1 <r_ 2$:
\begin{lem}[Basic estimate]\label{L:BE}
We have the following basic estimate:
$$\frac{M{_B}(0,r_1)}{\sqrt{K_B (0,r_2)}}= \frac{(n+1)\sqrt{\frac{n!}{\pi^n}\frac{1}{r_1 ^{n+2}}}}{\sqrt{\frac{n!}{\pi ^n}\frac{1}{r_2^{2n}}}}$$
$$ = \frac{(n+1)\sqrt{\frac{1}{r_1 ^{n+2}}}}{\frac{1}{r_2 ^n}} = (n+1)\sqrt{\frac{r_2 ^{2n}}{r_1 ^{n+2}}}\leq (n+1)\sqrt{\frac{(1+Cd)^{2n}}{(1-Cd)^{n+2}}}\leq (n+1)(1+\tilde d)\sqrt{\frac{|\xi_N|^2}{4d^2}+\frac{L(\xi_T)}{d}}$$
\end{lem}
Putting these estimates together and choosing the coordinates on the extremal domain
such that
$$p=0$$
is an extremal point such that
   the negative $x_1$-axis hits the boundary at 0 normally from the inside,
  then our estimates of the Bergman metric become:
$$\frac{M{_B}(0,r_2)}{\sqrt{K_B (0,r_1)}} = (n+1)\sqrt{\frac{r_1{{^{2n}}}}{r_2 ^{n+2}}} \geq (n+1) \sqrt{\frac{(1-Cd){{^{2n}}}}{(1+Cd){^{n+2}}}} $$
$$ \geq (n+1) (1-Cd)$$
for $\mathcal{C}^4$-boundaries; if the boundary is only $\mathcal{C}{^3}$, we only get the estimate
$$\geq (n+1) (1-C\sqrt{d})$$
Alltogether we have proved our main estimate for the Bergman metric:
\begin{thm}[Sharp estimates for the Bergman metric]
Put $B:= B{_\Omega} ((r,0),\xi)$, where $0\in b\Omega$ and $r>0$ is a small positive radius on the inner normal to $b\Omega$ at $0$. Then we have if the boundary is $\mathcal{C}^4,$
$$B\leq (n+1) (1+Cd) \sqrt{\frac{|\xi_N|^2}{4d^2}+\frac{L(\xi_T)}{d}}$$
In the opposite sense we have the estimates
\begin{equation}
B \geq
\begin{cases}
(n+1)(1-Cd)\sqrt{{\frac{|\xi_N|^2}{4d^2}+\frac{L(\xi_T)}{d}}} &\text{(if the boundary regularity is $\mathcal{C}^4 $})\\
\end{cases}
\end{equation}
\end{thm}

Klas Diederich\\
Bergische Universitat\\
Mathematik\\
Gausstr. 20 \\
42119 Wuppertal\\
Email: Klas.Diederich@math.uni-wuppertal.de\bg

J. E. Forn\ae ss\\
Department of Mathematical Sciences\\
Alfred Getz vei 1\\
Email: johnefo@math.ntnu.no

\end{document}